\documentclass[11pt]{article}
\usepackage[T2A]{fontenc}
\usepackage[cp1251]{inputenc}
\usepackage{amsfonts}
\usepackage{enumerate}
\usepackage{amssymb}

\usepackage{amsmath,amsthm}
\usepackage{amssymb,latexsym}
\usepackage{enumerate}

\begin{document}

\centerline{\textbf{On a group analogue of the Heyde theorem}}

\bigskip

\centerline{\textbf{Margaryta Myronyuk}}

\bigskip

\centerline{\textit{B. Verkin Institute for Low Temperature Physics
and Engineering }}

\centerline{\textit{of the National Academy of Sciences of Ukraine,
}}

\centerline{\textit{47 Nauky Ave, Kharkiv, 61103, Ukraine}}

\centerline{{Email: myronyuk@ilt.kharkov.ua}}

\begin{abstract}
Heyde proved that a Gaussian distribution on a real line is
characterized by the symmetry of the conditional distribution of one
linear form given another. The present article is devoted to an
analog of the Heyde theorem in the case when random variables take
values in a locally compact Abelian group and the coefficients of
the linear forms are integers.
\end{abstract}



\bigskip

\section{Introduction}

In 1953 V.P.Skitovich and G.Darmois proved independently that a
Gaussian distribution on a real line is characterized by the
independence of two linear forms of $n$ independent random variables
(\cite[$\S$3.1]{KaLiRa}). In 1970 C.C.Heyde proved a similar result
where a Gaussian distribution is characterized by the symmetry of
the conditional distribution of one linear form given another.

\bigskip

\textbf{The Heyde theorem} (\cite{He}).  \textit{Let $\xi_j, \ j=1,
2, ..., n, \ n \ge 2,$ be independent random variables. Consider
linear forms $L_1=a_1\xi_1+ \cdots +a_n\xi_n$ and $L_2=b_1\xi_1 +
\cdots + b_n\xi_n$, where the coefficients $a_j, b_j$ are nonzero
real numbers such that $b_i a_j + b_j a_i \ne 0$ for all $i, j$. If
the conditional distribution of $L_2$ given $L_1$ is symmetric then
all random variables $\xi_j$ are Gaussian.}

\bigskip

The Skitovich-Darmois theorem and the Heyde theorem were proved
using the finite difference method.

The Skitovich-Darmois theorem and the Heyde theorem were generalized
on locally compact Abelian groups (see e.g. \cite{Fe13},
\cite{Fe18}, \cite{FeGr2010}, \cite{MiFe1}, \cite{My2010},
\cite{My2013}). Specifically, in the article \cite{Fe7} G.Feldman
proved a generalization of the Skitovich-Darmois theorem in the case
when random variables take values in a locally compact Abelian group
and coefficients of the linear forms are integers. The main result
of this article is a generalization of the Heyde theorem in the case
when random variables take values in a locally compact Abelian group
and coefficients of the linear forms are integers (see \S3). To
obtain this result we use methods, which differ from methods of the
article \cite{Fe7}.

In the article  \cite{KS} A.M.Kagan and G.J.Sz\'ekely introduced a
notion of $Q$-independence of random variables which generalizes a
notion of independence of random variables. In particular, they
proved that some classical characterization theorems of mathematical
statistics hold true if instead of independence $Q$-independence is
considered. The article \cite{KS} has stimulated a series of
studies. Generalizations of some characterization theorems on
locally compact Abelian groups for $Q$-independent random variables
were obtained in \cite{Fe2017}, \cite{My2019}, \cite{My2018}. In \S
4 we prove that the results of \S3 hold true if instead of
independence $Q$-independence is considered.

To prove the main results of this article, we use the
finite-difference-method.

\section{Notation and definitions}

In the article we use standard results on abstract harmonic analysis
(see e.g. \cite{HeRo1}).

Let $X$ be a second countable locally compact  Abelian group,
$Y=X^\ast$ be its character group, and  $(x,y)$ be the value of a
character $y \in Y$ at an element $x \in X$. Let $H$ be a subgroup
of $Y$. Denote by $A(X,H)=\{x \in X: (x,y)=1 \ \ \forall \ y \in
H\}$ the annihilator of $H$. For each integer $n$, $n \ne 0,$ let
$f_n : X \mapsto X$ be the endomorphism $f_n x=nx.$ Set $X^{(n)} =
f_n(X)$, $X_{(n)}=Ker f_n$. Denote by ${\mathbb{Z}}(n)$  the finite
cyclic group of order $n$. For a fixed prime $p$ denote by
${\mathbb{Z}}(p^{\infty})$ the set of rational numbers of the form
$\{ k/p^n : k = 0, 1,..., p^n-1, n = 0, 1,...\}$ and define the
operation in ${\mathbb{Z}}(p^{\infty})$ as addition modulo 1. Then
${\mathbb{Z}}(p^{\infty})$ is transformed into an Abelian group,
which we consider in the discrete topology. Denote by $\Delta_p$ the
group of $p$-adic integers (\cite[\S 10.2]{HeRo1}). Note that
$\Delta_p^*\approx {\mathbb Z}(p^\infty)$ (\cite[\S 25.2]{HeRo1}).

Let ${M^1}(X)$ be the convolution semigroup of probability
distributions on $X$,  $\widehat \mu(y) = \int_X (x, y) d\mu(x)$ be
the characteristic function of a distribution $\mu \in {M^1}(X)$,
and $\sigma(\mu)$ be the support of $\mu$. If $H$ is a closed
subgroup of $Y$ and $\widehat \mu(y)=1$ for $y \in H$, then
$\widehat\mu(y+h) = \widehat\mu(y)$ for all $y \in Y, \ h \in H$ and
$\sigma(\mu) \subset A(X, H)$. For $\mu \in {M^1}(X)$ we define the
distribution $\bar \mu \in M^1(X)$ by the rule $\bar \mu(B) =
\mu(-B)$ for all Borel sets $B \subset X$. Observe that $\widehat
{\bar \mu}(y) = \overline{\widehat \mu(y)}$.

Let $x\in X$. Denote by $E_x$ the degenerate distribution
concentrated at the point $x$, and by $D(X)$ the set of all
degenerate distributions on
 $X$. A distribution $\gamma \in {M^1}(X)$ is called Gaussian (\cite[\S
4.6]{Pa}) if its characteristic function can be represented in the
form
$$ \widehat\gamma(y)= (x,y)\exp\{-\varphi(y)\},$$
where $x \in X$ and $\varphi(y)$ is a continuous nonnegative
function satisfying the equation
\begin{equation}\label{fe1}
    \varphi(u+v)+\varphi(u-v)=2[\varphi(u)+
\varphi(v)], \quad u, \ v \in Y.
\end{equation}
 Denote by $\Gamma(X)$ the set
of Gaussian distributions on $X$. We note that according to this
definition $D(X)\subset \Gamma(X)$. Denote by $I(X)$ the set of
shifts of Haar distributions $m_K$ of compact subgroups $K$ of the
group $X$. Note that

\begin{equation}\label{i1}
    \widehat{m}_K(y)=
\left\{%
\begin{array}{ll}
    1, & \hbox{$y\in A(Y,K)$;} \\
    0, & \hbox{$y\not\in A(Y,K)$.} \\
\end{array}%
\right.
\end{equation}
We note that if a distribution $\mu \in \Gamma(X)*I(X)$, i.e.
$\mu=\gamma*m_K$, where $\gamma\in \Gamma(X)$, then $\mu$ is
invariant with respect to the compact subgroup $K \subset X$ and
under the natural homomorphism $X \mapsto X/K \ \mu$ induces a
Gaussian distribution on the factor group $X/K$. Therefore the class
$\Gamma(X)*I(X)$ can be considered as a natural analogue of the
class $\Gamma(X)$ on locally compact Abelian groups.

Let $f(y)$ be a function on $Y$, and $h\in Y.$ Denote by $\Delta_h$
the finite difference operator
    $$\Delta_h f(y)=f(y+h)-f(y).$$
A function $f(y)$ on $Y$ is called a polynomial if
    $$\Delta_{h}^{n+1}f(y)=0$$
for some $n$ and for all $y,h \in Y$.

An integer $a$ is said to be admissible for a group $X$ if $X^{(a)}
\ne \{0\}$. The admissibility of integers $a_1, \dots , a_n$
 when we consider  the linear form  $a_1\xi_1 + \dots + a_n\xi_n$
 is a group analogue of the condition  $a_j \ne 0, \
j = 1, 2,\dots, n,$ for the case of  $X = \mathbb{R}$.

\section{The Heyde theorem for locally compact Abelian groups}

The main result of this section is  a full description of locally
compact Abelian  groups for which the group analogue of the Heyde
theorem takes place in the case when the characteristic functions of
considering random variables do not vanish.

\medskip

\textbf{Theorem 1.} \textit{Let $X$ be a second countable locally
compact Abelian group such that $X\neq X_{(2)}$. Then the following
statements hold:}

(I) \textit{Let $\xi_j, \ j=1, 2, ..., n, \ n \ge 2,$ be independent
random variables with values in $X$ and   distributions $\mu_j$ with
non-vanishing characteristic functions. Consider the linear forms
$L_1=a_1\xi_1+ \cdots +a_n\xi_n$ and $L_2=b_1\xi_1 + \cdots +
b_n\xi_n$, where the coefficients $a_j, \ b_j$ are admissible
integers for $X$ such that $b_i a_j+ b_j a_i$ are admissible
integers for $X$ for all $i, j$. Assume that the conditional
distribution of $L_2$ given $L_1$ is symmetric. Then the following
statements hold:}

\qquad \textit{$(i)$ If $X$ is a torsion-free group then all
$\mu_j\in\Gamma(X)$;}

\qquad \textit{$(ii)$ If $X=X_{(p)}$, where $p$ is a prime number
($p>2$), then all $\mu_j\in D(X)$.}

(II) \textit{If $X$ is not isomorphic to any of the groups mentioned
in (I), then there exist independent random variables $\xi_j, j=1,
2, ..., n, n \ge 2,$ with values in $X$ and distributions $\mu_j$
with non-vanishing characteristic functions, and admissible integers
$a_j, \ b_j $ such that $b_i a_j + b_j a_i$ are admissible integers
for $X$ for all $i, j$, such that the conditional distribution of
$L_2=b_1\xi_1 + \cdots + b_n\xi_n$ given $L_1=a_1\xi_1+ \cdots
+a_n\xi_n$ is symmetric, but all $\mu_j\not\in\Gamma(X)$.}

\medskip

\textbf{Remark 1}. Let $X$ be a locally compact Abelian group such
that every element of $X$ different from zero has order $p$, where
$p$ is a fixed prime number. Then $X$ is topologically isomorphic to
the group
\begin{equation}\label{rem1}
    \mathbb{Z}(p)^\mathfrak{m} \times \mathbb{Z}(p)^{\mathfrak{n}*},
\end{equation}
where $\mathfrak{m}$ and $\mathfrak{n}$ are arbitrary cardinal
numbers, ${\mathbb{Z}(p)}^{\mathfrak{m}}$ is considered in the
product topology, and $ {\mathbb{Z}(p)}^{\mathfrak{n}*}$ is
considered in the discrete topology (\cite[\S 25.29]{HeRo1}). If a
group $X$ is of the form (\ref{rem1}) for $p=2$ then there are no
exist admissible integers $a_j, \ b_j $ such that $b_i a_j + b_j
a_i$ are admissible integers for all $i, j$ for a group $X$. If a
group $X$ is of the form (\ref{rem1}) for $p=3$ then the set of
admissible integers $a_j, \ b_j $ such that $b_i a_j + b_j a_i$ are
admissible integers for all $i, j$ for a group $X$ is either
$a_j=b_j=1$ for all $j$ or $a_j=b_j=2$ for all $j$.

\medskip

To prove Theorem 1 we need some lemmas.

The following lemma was proved earlier for the case when the
coefficients of linear forms are topological automorphisms of the
group (\cite[\S 16.1]{FeBook2}). In the case when the coefficients
are integers, the proof is the same. For completeness of the
presentation, we give it in the article.

\medskip

\textbf{Lemma 1}. \textit{Let $\xi_j, \ j=1, 2, ..., n, \ n \ge 2,$
be independent random variables with values in a second countable
locally compact Abelian group $X$ and with distributions $\mu_j$.
Let $a_j, \ b_j$ be integers. The conditional distribution of the
linear form $L_2=b_1\xi_1 + \cdots + b_n\xi_n$ given $L_1=a_1\xi_1+
\cdots +a_n\xi_n$ is symmetric if and only if the characteristic
functions $\widehat{\mu}_j(y)$ satisfy the equation
\begin{equation}\label{1.1}
    \prod_{j = 1}^{n}{\widehat{\mu}_j(a_j u + b_j v)} =
    \prod_{j = 1}^{n}{\widehat{\mu}_j(a_j u - b_j v)},
\quad u,v\in Y.
\end{equation}
}

\medskip

\textbf{Proof.} Let $(\Omega, \mathfrak{A}, \mathrm{P})$ be a
probabilistic space, where the random variables $\xi_j$ are defined.
The condition of the symmetry of the conditional distribution of
$L_2$ given $L_1$ is equivalent to the equality
$${\rm P}\{\omega \in \Omega: L_1(\omega)\in A, L_2(\omega)\in B\}=
{\rm P}\{\omega \in \Omega: L_1(\omega)\in A, L_2(\omega)\in -B\} $$
for all Borel sets $A, B \subset X$. This means that the random
variables $(L_1,L_2)$ and $(L_1,-L_2)$ are identically distributed.
It is equivalent to the equality
$$
    \mathbf{E} \left[ (L_1, u ) (L_2,v) \right] =
    \mathbf{E} \left[ (L_1,u ) (-L_2,v) \right],
    \quad u,v\in Y.
$$

We obtain from the form of the linear forms that

$$  \mathbf{E} \left[ (a_1\xi_1+ \cdots +a_n, u ) (b_1\xi_1 + \cdots
+ b_n\xi_n,v) \right] =$$$$=
    \mathbf{E} \left[ (a_1\xi_1+ \cdots +a_n,u ) (-(b_1\xi_1 + \cdots + b_n\xi_n) \xi_j,v) \right],
\quad u,v\in Y, \Longleftrightarrow$$

$$  \mathbf{E} \left[ \prod_{j = 1}^{n}(a_j \xi_j, u ) \prod_{j = 1}^{n}(b_j \xi_j,v) \right] =
    \mathbf{E} \left[ \prod_{j = 1}^{n}(a_j \xi_j,u ) \prod_{j = 1}^{n}(-b_j \xi_j,v) \right],
\quad u,v\in Y, \Longleftrightarrow$$

$$  \mathbf{E} \left[ \prod_{j = 1}^{n}(\xi_j,a_j u ) \prod_{j = 1}^{n}(\xi_j,b_j v) \right] =
    \mathbf{E} \left[ \prod_{j = 1}^{n}(\xi_j,a_j u ) \prod_{j = 1}^{n}(\xi_j,-b_j v) \right],
\quad u,v\in Y, \Longleftrightarrow$$

\begin{equation}\label{l1.1}  \mathbf{E} \left[\prod_{j = 1}^{n}(\xi_j,a_j u + b_j v)\right] =
    \mathbf{E} \left[\prod_{j = 1}^{n}(\xi_j,a_j u - b_j v)\right],
\quad u,v\in Y.
\end{equation}

Taking into account the independence of the random variables
$\xi_j$, we have

\begin{equation}\label{l1.2}
    \prod_{j = 1}^{n}\mathbf{E}[(\xi_j,a_j u + b_j v)] =
    \prod_{j = 1}^{n}\mathbf{E}[(\xi_j,a_j u - b_j v)],
\quad u,v\in Y.
\end{equation}

Taking into account that $\widehat\mu_j(y)=\mathbf{E}[(\xi_j,y)]$,
we obtain (\ref{1.1}).
 $\blacksquare$

\medskip

The following lemma is an analogue of the Cramer theorem on the
decomposition of a Gaussian distribution on locally compact Abelian
groups.

\medskip

\textbf{Lemma 2} (\cite{Fe3}, see also \cite[\S 4.6]{FeBook2}).
\textit{Let $X$ be a second countable locally compact Abelian group.
Let $\mu=\mu_1 *\mu_2$, where $\mu\in \Gamma(X)$, $\mu_1, \mu_2 \in
M^1(X)$. If $X$ contains no subgroup topologically isomorphic to the
circle group $\mathbb{T}$ then $\mu_1, \mu_2\in \Gamma(X)$.}

\medskip

The following lemma is an analogue of the Marcinkiewicz theorem on
locally compact Abelian groups.

\medskip

\textbf{Lemma 3} (\cite{Fe6}, see also \cite[\S 5.11]{FeBook2}).
\textit{Let $X$ be a second countable locally compact Abelian group.
Let $\mu \in M^1(X)$ and its characteristic function is of the form
$$\widehat{\mu}(y)=\exp\{\psi(y)\},\quad \psi(0)=0,\quad y\in Y,$$ where $\psi(y)$
is a continuous polynomial. If $X$ contains no subgroup
topologically isomorphic to the circle group $\mathbb{T}$ then
$\mu\in \Gamma(X)$.}

\medskip

\textbf{Proof of Theorem 1.} {(I)} It follows from Lemma 1 that the
characteristic functions $\widehat\mu_j(y)$ satisfy equation
(\ref{1.1}).

Put $\nu_j=\mu_j*\bar\mu_j$. Then
$\widehat\nu_j(y)=|\widehat\mu_j(y)|^2> 0$. Obviously, the
characteristic functions $\widehat\nu_j(y)$ also satisfy equation
(\ref{1.1}). If we prove that $\nu_j\in \Gamma(X)$, then applying
Lemma 2, we obtain that $\mu_j\in \Gamma(X)$. Hence we can assume
from the beginning that all $\widehat\mu_j(y)>0$.

%Покажемо, що $\widehat\mu_j(y)$ --- характеристичні функції
%розподілів Гауса.

Set $\psi_j(y)=-\log \widehat\mu_j(y)$. We conclude from (\ref{1.1})
that

\begin{equation}\label{1.2}
    \sum_{j = 1}^{n}{\psi_j(a_j u + b_j
v)} = \sum_{j = 1}^{n}{\psi_j(a_j u - b_j v)}, \quad u,v\in Y.
\end{equation}

\noindent We use the finite difference method to solve equation
(\ref{1.2}).

Let $h_n$ be an arbitrary element of the group $Y$. Substitute
$u+b_n h_n$ for $u$ and $v+a_n h_n$ for $v$ in equation (\ref{1.2}).
Subtracting equation (\ref{1.2}) from the resulting equation we
obtain

\begin{equation}\label{1.3}
    \sum_{j = 1}^{n} \Delta_{l_{nj}h_n}
    \psi_j (a_j u + b_j v) = \sum_{j = 1}^{n-1} \Delta_{m_{nj}h_n}
    \psi_j(a_j u - b_j v), \quad u,v\in Y,
\end{equation}
where $l_{nj}=a_j b_n + b_j a_n$, $m_{nj}=a_j b_n - b_j a_n$. Note
that the right-hand side of equation (\ref{1.3}) does not contain
the function $\psi_n$.

Let $h_{n-1}$ be an arbitrary element of the group $Y$. Substitute
$u+b_{n-1} h_{n-1}$ for $u$ and $v+a_{n-1} h_{n-1}$ for $v$ in
equation (\ref{1.3}). Subtracting equation (\ref{1.3}) from the
resulting equation we obtain

\begin{equation}\label{1.4}
    \sum_{j = 1}^{n} \Delta_{l_{n-1,j}h_{n-1}} \Delta_{l_{nj}h_n}
    \psi_j (a_j u + b_j v) = \sum_{j = 1}^{n-2} \Delta_{m_{n-1,j}h_{n-1}} \Delta_{m_{nj}h_{n}}
    \psi_j(a_j u - b_j v), \quad u,v\in Y,
\end{equation}
where $l_{n-1,j}=a_j b_{n-1} + b_j a_{n-1}$, $m_{n-1,j}=a_j b_{n-1}
- b_j a_{n-1}$.

The right-hand side of equation (\ref{1.4}) does not contain the
functions $\psi_n$ and $\psi_{n-1}$. Arguing as above we get through
$n$ steps the equation

\begin{equation}\label{1.5}
    \sum_{j = 1}^{n} \Delta_{l_{1j} h_1 } \Delta_{l_{2j}h_2} ... \Delta_{l_{nj}h_n}
    \psi_j (a_j u + b_j v) = 0, \quad u,v\in Y,
\end{equation}
where $l_{ij}=a_j b_i + b_j a_i$.

Set $m_{ij}=a_j b_{i} - b_j a_{i}$ for all $i, j$.

Let $k_n$ be an arbitrary element of the group $Y$. Substitute
$u+b_n k_n$ for $u$ and $v-a_n k_n$ for $v$ in equation (\ref{1.5}).
Subtracting equation (\ref{1.5}) from the resulting equation we
obtain

\begin{equation}\label{1.6}
    \sum_{j = 1}^{n-1} \Delta_{m_{nj} k_n}
    \Delta_{l_{1j} h_1 } \Delta_{l_{2j}h_2} ... \Delta_{l_{nj}h_n}
    \psi_j (a_j u + b_j v) = 0, \quad u,v\in Y.
\end{equation}
Note that the left-hand side of equation (\ref{1.6}) does not
contain the function $\psi_n$. Arguing as above we get through $n-1$
steps the equation

\begin{equation}\label{1.7}
    \Delta_{m_{21} k_2} ... \Delta_{m_{n1} k_n}
    \Delta_{l_{11}h_1} \Delta_{l_{2,1}h_2} ... \Delta_{l_{n1}h_n}
    \psi_1 (a_1 u + b_1 v) = 0, \quad u,v\in Y.
\end{equation}

Put $v=0$ in (\ref{1.7}). We have

\begin{equation}\label{1.8}
    \Delta_{m_{21} k_2} ... \Delta_{m_{n1} k_n}
    \Delta_{l_{11}h_1} \Delta_{l_{2,1}h_2} ... \Delta_{l_{n1}h_n}
    \psi_1 (a_1 u ) = 0, \quad u\in Y.
\end{equation}

\medskip

We divide the proof into some steps. First, we prove the theorem,
assuming that all integers $b_i a_j - b_j a_i$ are admissible for
all $i \ne j$ for the group $X$. Then we prove that the case when
not all integers $b_i a_j - b_j a_i$ are admissible for all $i \ne
j$ for the group $X$ can be reduced to the case when all integers
$b_i a_j - b_j a_i$ are admissible for all $i \ne j$ for the group
$X$.

\medskip

\textbf{1.} Let all integers $b_i a_j - b_j a_i$ be admissible for
all $i \ne j$ for the group $X$.

(i) Let $X$ be a torsion-free group, i.e. $X_{(m)}=\{0\}$ for all
$m\in \mathbb{Z}$, $m\ne 0$. This implies that
$\overline{Y^{(m)}}=Y$ for all integers $m, \ m\ne 0$
(\cite[24.41]{HeRo1}). Particularly, $\overline{Y^{(a_j)}}=Y$, $j=1,
2, ..., n$. Since integers $l_{ij}, m_{ij}$ are admissible, it means
for torsion-free groups that $l_{ij}\neq 0, m_{ij} \neq 0$.
Therefore, $\overline{Y^{(l_{ij})}}=Y$, $\overline{Y^{(m_{ij})}}=Y$,
$j=1, 2, ..., n$. Taking into account that $h_j, k_j$ are arbitrary
elements of the group $Y$, (\ref{1.8}) yields that the function
$\psi_1(y)$ satisfy the equation

\begin{equation}\label{1.9}
    \Delta_{h}^{2n-1}
    \psi_1 (y) = 0, \quad y, h\in Y.
\end{equation}

So $\psi_1 (y)$ is a continuous polynomial on the group $Y$. It
follows from Lemma 3 that $\mu_1\in\Gamma(X)$. Arguing as above we
prove that $\mu_j\in\Gamma(X)$ ($j=2,\dots,n$).

(ii) Let $ X=X_{(p)} $, where $p$ is a prime number ($p\neq 2$).

Since integers $a_j, b_j$, $j=1, ... n$, $b_i a_j \pm b_j a_i$ are
admissible for the group $X$, we get that endomorphisms $f_{a_j},
f_{b_j}, f_{b_i a_j \pm b_j a_i}$ are automorphisms. It means that
$Y^{(a_j)}=Y$, $Y^{(b_j)}=Y$, $Y^{(b_i a_j \pm b_j a_i)}=Y$. Taking
into account that $h_j, k_j$ are arbitrary elements of $Y$, it
follows from the form of integers $l_{ij}, m_{ij}$ and (\ref{1.8})
that the function $\psi_1(y)$ satisfy equation (\ref{1.9}). Arguing
as in case 1(i), we get that all $\mu_j\in\Gamma(X)$. Since that
component of zero of the group $X$ is equal to zero, we have in case
(ii) that $\Gamma(X)=D(X)$.

\medskip

\textbf{2.} Suppose now that for some $i,j$ integers $b_i a_j - b_j
a_i$ are admissible for $X$ and for some $i,j$ integers $b_i a_j -
b_j a_i$ are not admissible for $X$.

%Припустимо, що $n=2$, тобто $L_1=a_1\xi_1+a_2\xi_2$ и $L_2=b_1\xi_1
%+ b_2\xi_2$.

(i) Let $X$ be a torsion-free group. In this case if integers $b_i
a_j - b_j a_i$ are not admissible for $X$ then $b_i a_j - b_j a_i=0$
for some $i,j$.

Renumbering random variables, we can assume that

\begin{equation}\label{1.9.1}
    {a_1 \over b_1}= \dots = {a_{r_1} \over b_{r_1}}=c_1, {a_{r_1+1} \over b_{r_1+1}}= \dots = {a_{r_2} \over b_{r_2}}=c_2,\dots,
    {a_{r_{k}+1} \over b_{r_{k}+1}}= \dots = {a_n \over b_n}=c_{k+1},
\end{equation}
where $1\leq r_1<r_2<\dots < r_k< n$, and $c_i\neq c_j$ for all $i
\ne j$.

Consider integers

\begin{equation}\label{1.9.2}
    d_0={b_1 b_{r_1+1}\dots b_{r_k+1}\over b_1}, \quad d_j={b_1
    b_{r_1+1}\dots b_{r_k+1}\over b_{r_j+1}}, \quad j=1,2,\dots, k,
\end{equation}
and random variables

$$\zeta_1=b_1\xi_1+\dots+b_{r_1}\xi_{r_1},\quad \zeta_2=b_{r_1+1}\xi_{r_1+1}+\dots+b_{r_2}\xi_{r_2}, \dots,
\quad \zeta_{k+1}=b_{r_k+1}\xi_{r_k+1}+\dots+b_{n}\xi_{n}.$$ Put
$$L_1'=a_1 d_0 \zeta_1 + a_{r_1+1} d_1 \zeta_2 + \cdots +a_{r_{k}+1}
d_k \zeta_{k+1}, \quad L_2'=\zeta_1 + \zeta_2 + \cdots
+\zeta_{k+1}.$$

Since $L_1'= b_1 b_{r_1+1}\dots b_{r_k+1} L_1$, $L_2'=L_2$, it is
obvious that the conditional distribution of $L_2'$ given $L_1'$ is
also symmetric. It follows from (\ref{1.9.1}) and (\ref{1.9.2}) that
the coefficients of the linear forms $L_1'$ and $L_2'$ satisfy
conditions $a_1 d_0 \neq 0, a_{r_1+1} d_1 \neq 0, \cdots,
a_{r_{k}+1} d_k \neq 0 $, $a_1 d_0 \pm a_{r_j+1} d_j \neq 0,
a_{r_i+1} d_i \pm a_{r_j+1} d_j \neq 0$. It follows from case 1(i)
that all $\zeta_j$ have Gaussian distributions. Applying Lemma 2 we
get that $b_1\xi_1, \cdots, b_n\xi_n$ have Gaussian distributions.
Since $X$ is a torsion-free group, the endomorphisms $f_{b_1},
\cdots, f_{b_n}$ are monomorphisms. Now it is easy to verify that
the random variables $\xi_1, \cdots, \xi_n$ have Gaussian
distributions.

(ii) Let $X=X_{(p)}$, where $p$ is a prime number ($p>2$). Since
integers $a_j, b_j$ are admissible for $X$, the endomorphisms
$f_{a_1}, \cdots, f_{a_n}, f_{b_1}, \cdots, f_{b_n}$ are
automorphisms. Therefore we can take into consideration new random
variables and assume that all $a_j=1$ and all $b_j \in \{ 1,2,\dots,
p-1 \} $. Then we have $b_i - b_j =0$ for some $i, j$. Arguing as in
case 2(i), we get that all $\xi_1, \cdots, \xi_n$ are degenerated.

\medskip

\textbf{3.} Let all integers $b_i a_j - b_j a_i$ be not admissible
for $X$. Then $b_1L_1=a_1L_2$. It is obvious that all random
variables $\xi_1, \cdots, \xi_n$ are degenerated.

\medskip

{(II)} Suppose now that $X$ is not isomorphic to any of the groups
mentioned in (I). Then $X$ contains an element $x_0$ of order $p$,
where $p$ is a prime number and $X^{(p)}\ne\{0\}$. Let $M$ be a
subgroup of $X$ generated by the element $x_0$. Then
$M\cong\mathbb{Z}(p)$. Let $\xi_1$ and $\xi_2$ be independent
identically distributed random variables with values in $M$ and with
a distribution

\begin{equation}\label{1.10.1}
    \mu= aE_0+(1-a)m_M,
\end{equation}
where $0<a<1$. We consider the distribution $\mu$ as a distribution
on the group $X$. Taking into account (\ref{i1}), we have that

\begin{equation}\label{1.10.2}
    \widehat{\mu}(y)= \left\{%
\begin{array}{ll}
    1, & \hbox{$y\in A(Y,M)$,} \\
    a, & \hbox{$y\not\in A(Y,M)$.} \\
\end{array}%
\right.
\end{equation}
We consider the linear forms $L_1=p\xi_1-\xi_2$ and
$L_2=\xi_1+p\xi_2$. Since $X^{(p)}\neq \{0\}$ and $X_{(p)}\neq
\{0\}$, the integers $p, p^2-1$ are admissible for $X$. By Lemma 1
the conditional distribution of $L_2$ given $L_1$ is symmetric if
and only if the characteristic function $\widehat\mu(y)$ satisfy
equation (\ref{1.1}) which takes the form

\begin{equation}\label{1.10}
    \widehat\mu(pu+v)\widehat\mu(-u+pv)=\widehat\mu(pu-v)\widehat\mu(-u-pv),
    \quad u, v \in Y.
\end{equation}
Since $\sigma(\mu) \subset M$, we have $\widehat\mu(y)=1$ for $y\in
A(Y,M)$ and $\widehat\mu(y)$ is $A(Y,M)$-invariant, i.e.
$\widehat\mu(y)$ takes a constant value on each coset of the group
$Y$ with respect to $A(Y,M)$. Then equation (\ref{1.10}) induces an
equation on the factor-group

\begin{equation}\label{1.11}
    \widehat\mu(p[u]+[v])\widehat\mu(-[u]+p[v])=\widehat\mu(p[u]-[v])\widehat\mu(-[u]-p[v]),
    \quad [u], [v] \in Y/A(Y,M).
\end{equation}
Since $Y/A(Y,M)\approx \mathbb{Z}(p)$, equation (\ref{1.11}) is
transformed into the equality. Hence by Lemma 1  the conditional
distribution of $L_2$ given $L_1$ is symmetric.
 $\blacksquare$

\medskip

Now we describe locally compact Abelian groups for which the group
analogue of the Heyde theorem takes place without assumption that
the characteristic functions of the considering random variables  do
not vanish.  Note that the obtained class of groups is changed
(compare with Theorem 1).

\medskip

\textbf{Theorem 2.} \textit{Let $X$ be a second countable locally
compact Abelian group such that $X\neq X_{(2)}$. Then the following
statements hold:}

(I) \textit{Let $\xi_j, \ j=1, 2, ..., n, \ n \ge 2,$ be independent
random variables with values in $X$ and   distributions $\mu_j$.
Consider the linear forms $L_1=a_1\xi_1+ \cdots +a_n\xi_n$ and
$L_2=b_1\xi_1 + \cdots + b_n\xi_n$, where the coefficients $a_j, \
b_j$ are admissible integers for $X$ such that $b_i a_j+ b_j a_i$
are admissible integers for $X$ for all $i, j$. Assume that the
conditional distribution of $L_2$ given $L_1$ is symmetric. Then the
following statements hold:}

\qquad \textit{$(i)$ If $X=\mathbb{R}^n\times D$, where $n\geq 0$
and $D$ is a torsion-free group, then all $\mu_j\in\Gamma(X)$;}

\qquad \textit{$(ii)$ If $X=X_{(3)}$ then all $\mu_j\in D(X)$.}

(II) \textit{If $X$ is not isomorphic to any of the groups mentioned
in (I), then there exist independent random variables $\xi_j, j=1,
2, ..., n, n \ge 2,$ with values in $X$ and distributions $\mu_j$,
and admissible integers $a_j, \ b_j $ such that $b_i a_j + b_j a_i$
are admissible integers for $X$ for all $i, j$, such that the
conditional distribution of $L_2=b_1\xi_1 + \cdots + b_n\xi_n$ given
$L_1=a_1\xi_1+ \cdots +a_n\xi_n$ is symmetric, but all
$\mu_j\not\in\Gamma(X)*I(X)$.}

\medskip

\textbf{Remark 2.} In Theorem 2 we described all groups on which the
symmetry of the conditional distribution of one linear form given
another implies that all distributions belong to the class
$\Gamma(X)*I(X)$. In fact if a distribution belongs to the class
$\Gamma(X)*I(X)$ then it belongs to the class $\Gamma(X)$.

\medskip

To prove Theorem 2 we need some lemmas.

\medskip

\textbf{Lemma 4.} \textit{Let $K$ be a connected compact Abelian
group and $\mu_j$ are distributions on the group $K^*$. Assume that
the characteristic functions $\widehat{\mu}_j(y)$  satisfy equation
$(\ref{1.1})$ on the group $K$  and $\widehat{\mu}_j(y)\geq 0$. If
$a_j b_j\neq 0 $, then all $\widehat\mu_j(y)=1$, $y\in K$.}

\medskip

The proof of Lemma 4 is carried out in the same way as the proof of
Lemma 1 of the article \cite{Fe7}, but it is based on the Heyde
theorem. For completeness of the statement we will give this proof
in the article.

\bigskip

\textbf{Proof of Lemma 4.} Two cases are possible.

1. $K\not\approx \mathbb{T}$. Since $K$ is a connected compact
group, there exists a continuous monomorphism $f:
\mathbb{R}\rightarrow K$ such that $\overline{f(\mathbb{R})}=K$
(\cite[\S 25.18]{HeRo1}). Consider the restriction of equation
(\ref{1.1}) to $f(\mathbb{R})$. It follows from the Heyde theorem
that $\widehat\mu_j(f(t))=\exp\{ -\lambda_j t^2\}$, where
$\lambda_j\geq 0$ for $t\in \mathbb{R}$.

Let $U$ be a neighborhood of zero of $K$. Since $f$ is a
monomorphism and $\overline{f(\mathbb{R})}=K$, there exists a
sequence $t_n\longrightarrow +\infty$ such that $f(t_n)\in U$ for
all $n$. If $\lambda_j>0$ for some $j$ then
$\widehat\mu_j(f(t))=\exp\{ -\lambda_j t^2 \}\longrightarrow 0$ for
$t_n\longrightarrow +\infty$, which contradicts to the continuity of
the function $\widehat\mu_j(y)$ at zero. We have that all
$\lambda_j= 0$. So, $\widehat\mu_j(f(t))=1$ for $t\in \mathbb{R}$.
Since $\overline{f(\mathbb{R})}=K$, we have $\widehat\mu_j(y)=1$ for
$y\in K$.

2. $K\approx \mathbb{T}$. In this case the characteristic functions
$\widehat{\mu}_j(y)$ are $2\pi$-periodic and satisfy equation
(\ref{1.1}). It follows from the Heyde theorem that
$\widehat\mu_j(y)$ are the characteristic functions of Gaussian
distributions, i.e. $\widehat\mu_j(y)=\exp\{ -\lambda_j y^2\}$,
where $\lambda_j\geq 0$. It follows from $2\pi$-periodicity that
$\lambda_j= 0$. So, $\widehat\mu_j(y)=1$ for $y\in K$.
$\blacksquare$

\medskip

\textbf{Lemma 5.} \textit{Let $X=\Delta_p$. Then there exist
independent random variables $\xi_1, \xi_2$ with values in $X$ and
distributions $\mu_1, \mu_2$ such that the conditional distribution
of $L_2=\xi_1 + p\xi_2$ given $L_1=p\xi_1- \xi_2$ is symmetric and
$\mu_1, \mu_2\not\in I(X)$.}

\medskip

\textbf{Proof.} Since $X=\Delta_p$, we have
$Y\approx\mathbb{Z}(p^{\infty})$. We can assume without loss of
generality that $Y=\mathbb{Z}(p^{\infty})$. Consider the
distribution (\ref{1.10.1}) on $\mathbb{Z}(p)$.   Note that
$Y_{(p)}\approx \mathbb{Z}(p)$. Obviously,

$$g(y)=\left\{%
\begin{array}{ll}
    1, & \hbox{$y=0$;} \\
    a, & \hbox{$y\ne 0$,} \\
\end{array}%
\right.     $$ is a characteristic function on the group $Y_{(p)}$.
Consider the function
$$f(y)=\left\{%
\begin{array}{ll}
    1, & \hbox{$y=0$;} \\
    a, & \hbox{$y\in Y_{(p)}\setminus \{0\}$;} \\
    0, & \hbox{$y\not\in Y_{(p)}$,} \\
\end{array}%
\right.    $$ where $0<a<1$. It follows from \cite[\S 32.43]{HeRo2}
that the function $f(y)$ is positive definite. By the Bochner
theorem (\cite[\S 33.3]{HeRo2}) the function $f(y)$ is the
characteristic function of a distribution $\nu$ on $X$. Put
$\mu_1=\mu_2=\nu$ and verify that the conditional distribution of
$L_2=\xi_1 + p\xi_2$ given $L_1=p\xi_1- \xi_2$ is symmetric. By
Lemma 1 it suffices to show that the characteristic function $f(y)$
satisfies equation (\ref{1.10}). Obviously, if $u,v \in Y_{(p)}$
then equation (\ref{1.10}) is satisfied. If either $u\in Y_{(p)}$,
$v \not \in Y_{(p)}$, or $u\not \in Y_{(p)}$, $v \in Y_{(p)}$, both
sides of equation (\ref{1.10}) are equal to zero. Let $u,v \not\in
Y_{(p)}$. If the left-hand side of equation (\ref{1.10}) is not
equal to zero, then $pu+v, -u+pv\in Y_{(p)}$. We conclude from this
that $(p^2+1)v\in Y_{(p)}$. Hence $v\in Y_{(p)}$ contrary to the
assumption.  Thus, the left-hand side of equation (\ref{1.10}) is
equal to zero. Arguing in the same way, we get that the right-hand
side of equation (\ref{1.10}) is equal to zero. We get that the
function $f(y)$ satisfies equation (\ref{1.10}). It is clear that
$\nu$ can be chosen such that $\mu\not\in I(X)$. $\blacksquare$

\medskip

\textbf{Remark 3.} If in Lemma 5 independent random variables have
distributions with non-vanishing characteristic functions then
Theorem 1 implies that the distributions are degenerated.

\medskip

\textbf{Lemma 6.} \textit{Let $X=\mathbb{Z}(p)$, where $p$ is a
prime number ($p>3$). There exist independent random variables
$\xi_1, \xi_2, \xi_3, \xi_4$ with values in $X$ and distributions
$\mu_1, \mu_2, \mu_3, \mu_4$ such that the conditional distribution
of $L_2=\xi_1 + \xi_2+2\xi_3 + 2\xi_4$ given $L_1=\xi_1 +
\xi_2+\xi_3 + \xi_4$ is symmetric and $\mu_1, \mu_2, \mu_3, \mu_4
\not\in I(X)$.}

\medskip

\textbf{Proof.} Since $X=\mathbb{Z}(p)$, we have
$Y\approx\mathbb{Z}(p)$. Choose nonzero elements $y_1, y_2\in Y$such
that $y_1\neq \pm y_2$. Let $\nu_1, \nu_2$ be distributions in $X$
with the densities $\rho_1(x)=1+{1\over 2} {\rm Re}(x,y_1)$,
$\rho_2(x)=1+{1\over 2} {\rm Re}(x,y_2)$ with respect to the Haar
distribution $m_X$ respectively. Then

\begin{equation}\label{6.1}
    \widehat{\nu}_1(y)= \left\{%
\begin{array}{ll}
    1, & \hbox{$y=0$;} \\
    {1\over 2}, & \hbox{$y=\pm y_1$;} \\
    0, & \hbox{$y\not \in \{0,\pm y_1\}$.} \\
\end{array}%
\right.
    \quad
    \widehat{\nu}_2(y)= \left\{%
\begin{array}{ll}
    1, & \hbox{$y=0$;} \\
    {1\over 2}, & \hbox{$y=\pm y_2$;} \\
    0, & \hbox{$y\not \in \{0,\pm y_2\}$.} \\
\end{array}%
\right.
\end{equation}
Put $\mu_1=\mu_3=\nu_1$, $\mu_2=\mu_4=\nu_2$ and verify that the
conditional distribution of $L_2=\xi_1 + \xi_2+2\xi_3 + 2\xi_4$
given $L_1=\xi_1 + \xi_2+\xi_3 + \xi_4$ is symmetric. By Lemma 1 it
suffices to show that the characteristic functions of distributions
$\mu_1, \mu_2, \mu_3, \mu_4$ satisfy equation (\ref{1.1}), which
takes the form

\begin{equation}\label{6.2}
    \widehat{\nu}_1(u+v)\widehat{\nu}_2(u+v)\widehat{\nu}_1(u+2v)\widehat{\nu}_2(u+2v)=
    \widehat{\nu}_1(u-v)\widehat{\nu}_2(u-v)\widehat{\nu}_1(u-2v)\widehat{\nu}_2(u-2v),
    \quad u,v \in Y.
\end{equation}
Since $y_1\neq \pm y_2$, taking into account (\ref{i1}) we have
$\widehat{\nu}_1(y)\widehat{\nu}_2(y)=\widehat{m}_X(y)$. Equation
(\ref{6.2}) is equivalent to the equation
\begin{equation}\label{6.3}
    \widehat{m}_X(u+v)\widehat{m}_X(u+2v)=
    \widehat{m}_X(u-v)\widehat{m}_X(u-2v),
    \quad u,v \in Y.
\end{equation}
Obviously, equation (\ref{6.3}) is satisfied if either $u=0$ or
$v=0$. Let $u\neq 0$, $v\neq 0$. If the left-hand side of equation
(\ref{6.3}) is not equal to zero, then $u+v=0$, $u+2v=0$. We
conclude from this that $v=0$ contrary to the assumption. Thus, the
left-hand side of equation (\ref{6.3}) is equal to zero. Arguing in
the same way, we get that the right-hand side of equation
(\ref{6.3}) is equal to zero. We get that the function
$\widehat{m}_X(y)$ satisfies equation (\ref{6.3}). It is clear that
$\nu_1, \nu_2\not\in I(X)$. $\blacksquare$

\bigskip

\textbf{Proof of Theorem 2.} {(I)}(i) If $X= \mathbb{R}^n\times D$,
then $Y\approx \mathbb{R}^n\times K$, where $K$ is a connected
compact group.  We can assume without loss of generality that
$Y=\mathbb{R}^n\times K$.

Lemma 1 implies that the characteristic functions of distributions
$\mu_j$ satisfy equation (\ref{1.1}). Reasoning as in beginning of
the proof of Theorem 1, we can assume that $\widehat\mu_j(y)\geq0$.

Consider the restriction of equation (\ref{1.1}) to $K$. It follows
from Lemma 4 that $\widehat\mu_j(y)=1$ on $K$. Then
$\sigma(\mu_j)\subset A(X,K)= \mathbb{R}^n$. We consider the
restriction of equation (\ref{1.1}) to each one-dimensional subspace
and obtain from the Heyde theorem that all restriction of the
characteristic functions $\widehat\mu_j(y)$ are the characteristic
functions of Gaussian distributions. It follows from this that
$\mu_j\in\Gamma(X)$.

(ii) If $X=X_{(3)}$ then we can assume without loss of generality
that $a_j, b_j$ are equal to $\pm 1$. Since all integers $b_i a_j +
b_j a_i$ are admissible, we have either $L_1=L_2$ or $L_1=-L_2$. It
easily follows from this that all random variables $\xi_1, \cdots,
\xi_n$ are degenerated.

{(II)} 1.  By the structure theorem for locally compact Abelian
groups $X\approx \mathbb{R}^n\times G$, where $G$ contains an open
compact subgroup $H$ (\cite[24.30]{HeRo1}). Suppose that $X$ is a
torsion-free group. Since $X$ is not as in case (I), we have $H\neq
\{0\}$. A compact torsion-free group is topologically isomorphic to
the group
$$H\approx (\Sigma_a)^\mathfrak{n}\times \sum_{p\in \mathcal{P}}
\Delta_p^{\mathfrak{n}_p},$$ where $\Sigma_a^*=\mathbb{Q}$,
$\mathfrak{n}, \mathfrak{n}_p$ are cardinal numbers (see
\cite[25.8]{HeRo1}). Note that for each prime number $p$ a group
$\Sigma_a$ contains a subgroup $\Delta_p$. It follows from the fact
that the factor-group $\mathbb{Q}/ \mathbb{Z}$ contains a subgroup
$\mathbb{Z}(p^{\infty})$ for each prime number $p$. Then the
statement of Theorem 2 follows from Lemma 5.

2. Suppose that $X$ contains an element $x_0$ of order $p$, where
$p$ is a prime number and $ X\ne X_{(p)}$. Then the statement of
Theorem 2 follows part (II) of Theorem 1.

3. Suppose that $ X= X_{(p)}$, where $p>3$. Then the statement of
Theorem 2 follows from Lemma 6.

Theorem 2 is proved.
 $\blacksquare$

\section{The Heyde theorem for locally compact Abelian groups for Q-independent
random variables}

In the article  \cite{KS} A.M.Kagan and G.J.Sz\'ekely introduced a
notion of $Q$-independence of random variables which generalizes the
notion of independence of random variables. Then in \cite{Fe2017}
G.M. Feldman in a natural way introduced a notion of
$Q$-independence of random variables taking values in a locally
compact Abelian group. He proved that if we consider
$Q$-independence instead of independence, then the group analogue of
the Cram\'er theorem about decomposition of a Gaussian distribution
and some group analogues of the the Skitovich--Darmois and Heyde
theorems hold true for the same classes of groups. These studies
were continued in \cite{My2019}. Some analogues of characterization
theorems for Q-independent random variables with values in Banach
spaces were studied in \cite{My2018}. We prove in this section that
the results of \S 3 hold true if instead of independence
$Q$-independence is considered.

Let $\xi_1, \dots, \xi_n$ be random variables with values in the
group $X$. We say that the random variables $\xi_1, \dots, \xi_n$
are $Q$-independent if the characteristic function of the vector
$(\xi_1, \dots, \xi_n)$ can be represented in the form
\begin{equation}\label{i0}
    \hat\mu_{(\xi_1, \dots, \xi_n)}(y_1, \dots, y_n)={\bf E}[(\xi_1,
y_1)\cdots(\xi_n,
y_n)]=$$$$=\left(\prod_{j=1}^n\hat\mu_{\xi_j}(y_j)\right)\exp\{q(y_1,
\dots, y_n)\}, \quad y_j\in Y,
\end{equation}
where $q(y_1, \dots, y_n)$ is a continuous polynomial on the group
$Y^n$ such that $q(0, \dots, 0)=0$.

As is well known, any continuous polynomial is equal to the constant
on compact elements (see e.g. \cite[\S5.7]{FeBook2}). If the group
$Y$ consists only of compact elements then the connected component
of zero of the group $ X $ is equal to zero (see e.g.
\cite[\S24.17]{HeRo1}). Thus, independence and Q-independence of
random variables are equivalent on groups $X$ whose connected
component of zero is equal to zero.

We prove that Theorem 1 remains true if we change the condition of
independence of $\xi_1, \dots, \xi_n$ for $Q$-independence. The
following statement holds true.

\medskip

\textbf{Theorem 3.} \textit{Let $X$ be a second countable locally
compact Abelian group such that $X\neq X_{(2)}$. Then the following
statements hold:}

(I) \textit{Let $\xi_j, \ j=1, 2, ..., n, \ n \ge 2,$ be
Q-independent random variables with values in $X$ and with
distributions $\mu_j$ with non-vanishing characteristic functions.
Consider the linear forms $L_1=a_1\xi_1+ \cdots +a_n\xi_n$ and
$L_2=b_1\xi_1 + \cdots + b_n\xi_n$, where the coefficients $a_j, \
b_j$ are admissible integers for $X$ such that $b_i a_j+ b_j a_i$
are admissible integers for $X$ for all $i, j$. Assume that the
conditional distribution of $L_2$ given $L_1$ is symmetric. Then the
following statements hold:}

\qquad \textit{$(i)$ If $X$ is a torsion-free group then all
$\mu_j\in\Gamma(X)$;}

\qquad \textit{$(ii)$ If $X=X_{(p)}$, where $p$ is a prime number
($p>2$), then all $\mu_j\in D(X)$.}

(II) \textit{If $X$ is not isomorphic to any of the groups mentioned
in (I), then there exist Q-independent random variables $\xi_j, j=1,
2, ..., n, n \ge 2,$ with values in $X$ and distributions $\mu_j$
with non-vanishing characteristic functions, and admissible integers
$a_j, \ b_j $ such that $b_i a_j + b_j a_i$ are admissible integers
for $X$ for all $i, j$ such that the conditional distribution of
$L_2=b_1\xi_1 + \cdots + b_n\xi_n$ given $L_1=a_1\xi_1+ \cdots
+a_n\xi_n$ is symmetric, but all $\mu_j\not\in\Gamma(X)$.}

\medskip

To prove Theorem 3 we need the following lemma.

\medskip

\textbf{Lemma 7}. \textit{Let $\xi_j, \ j=1, 2, ..., n, \ n \ge 2,$
be Q-independent random variables with values in a second countable
locally compact Abelian group $X$ and with distributions $\mu_j$.
Let $a_j, \ b_j$ be integers. The conditional distribution of the
linear form $L_2=b_1\xi_1 + \cdots + b_n\xi_n$ given $L_1=a_1\xi_1+
\cdots +a_n\xi_n$ is symmetric if and only if the characteristic
functions $\widehat{\mu}_j(y)$ satisfy the equation
\begin{equation}\label{1.1q}
    \prod_{j = 1}^{n}{\widehat\mu_j(a_j u + b_j v)} =
    \prod_{j = 1}^{n}{\widehat\mu_j(a_j u - b_j v)} \exp\{q(u,v)\},
\quad u,v\in Y,
\end{equation}
where $q(u,v)$ is a continuous polynomial on the group $Y^2$,
$q(0,0)=0$.}

\medskip

\textbf{Proof.} The condition of the symmetry of the conditional
distribution of the linear form $L_2=b_1\xi_1 + \cdots + b_n\xi_n$
given $L_1=a_1\xi_1+ \cdots +a_n\xi_n$ is equivalent to equation
(\ref{l1.1}).

Taking into account the Q-independence of the random variables
$\xi_j$, we have

\begin{equation}\label{l1.2}
    \prod_{j = 1}^{n}\mathbf{E}[(\xi_j,a_j u + b_j v)] \exp\{q_0(a_1 u + b_1 v,\dots, a_n u + b_n v)\}
    =$$$$=
    \prod_{j = 1}^{n}\mathbf{E}[(\xi_j,a_j u - b_j v)] \exp\{q_0(a_1 u - b_1 v,\dots, a_n u - b_n v)\},
\quad u,v\in Y,
\end{equation}
where $q_0(y_1,\dots, y_n)$  is a continuous polynomial on the group
$Y^n$, $q_0(0,\dots, 0)=0$. Put $q(u,v)=q_0(a_1 u - b_1 v,\dots, a_n
u - b_n v)-q_0(a_1 u + b_1 v,\dots, a_n u + b_n v)$. Taking into
account that $\widehat\mu_j(y)=\mathbf{E}[(\xi_j,y)]$, we obtain
equation (\ref{1.1q}).
 $\blacksquare$

\medskip

\textbf{Proof of Theorem 3.} {(I)} Lemma 7 implies that the
characteristic functions of distributions $\mu_j$ satisfy equation
(\ref{1.1q}). Reasoning as in the beginning of the proof of Theorem
1, we can assume that $\widehat\mu_j(y)>0$ and the polynomial
$q(u,v)$ is real-valued.

We will show that $\widehat\mu_j(y)$ are the characteristic
functions of Gaussian distributions. Put $\psi_j(y)=-\log
\widehat\mu_j(y)$. It follows from (\ref{1.1q}) that

\begin{equation}\label{1.2q}
    \sum_{j = 1}^{n}{\psi_j(a_j u + b_j
v)} = \sum_{j = 1}^{n}{\psi_j(a_j u - b_j v)} + q(u,v), \quad u,v\in
Y.
\end{equation}

\noindent We use the finite-difference method in the same manner as
in Theorem 1 to solve equation (\ref{1.2q}). We will retain the
notation from the proof of Theorem 1. In fact we add one more step.
Instead of equation (\ref{1.7}), we get through $2n-1$ steps the
equation

\begin{equation}\label{5.7q}
    \Delta_{m_{11} k_2} ... \Delta_{m_{n1} k_n}
    \Delta_{l_1,1} \Delta_{l_2,1} ... \Delta_{l_{n-1,1}} \Delta_{l_{n1}}
    \psi_1 (a_1 u + b_1 v) =$$$$ =
    \Delta_{(b_2 h_2, -a_2 h_2)}...\Delta_{(b_n h_n, -a_n h_n)}\Delta_{(b_1 h_1, a_1 h_1)}...
    \Delta_{(b_n h_n, a_n h_n)} q(u,v), \quad u,v\in Y.
\end{equation}

Since $q(u,u)$ is a polynomial, we have
\begin{equation}\label{5.7.1}
  \Delta^{l+1}_{(h, k)}q(u,v)=0,
\quad u, v\in Y,
\end{equation}
for some $l$ and arbitrary elements $h$ and $k$ of $Y$.

We substitute $v=0$ in (\ref{5.7q}) and apply the operator
$\Delta^{l+1}_{(h, k)}$ to (\ref{5.7q}). Taking into account
(\ref{5.7.1}), we obtain

\begin{equation}\label{5.8q}
    \Delta^{l+1}_{h}\Delta_{m_{11} k_2} ... \Delta_{m_{n1} k_n}
    \Delta_{l_1,1} \Delta_{l_2,1} ... \Delta_{l_{n-1,1}} \Delta_{l_{n1}}
    \psi_1 (a_1 u ) = 0, \quad u\in Y.
\end{equation}

We complete the proof of Theorem 2 in the same way as Theorem 1.

(II) Since independent random variables are Q-independent, this part
of Theorem 3 follows from part (II) of Theorem 1. $\blacksquare$

\bigskip

We prove that Theorem 2 remains true if we change the condition of
independence of $\xi_1, \dots, \xi_n$ for $Q$-independence. The
following statement holds true.

\medskip

\textbf{Theorem 4.} \textit{Let $X$ be a second countable locally
compact Abelian group such that $X\neq X_{(2)}$. Then the following
statements hold:}

(I) \textit{Let $\xi_j, \ j=1, 2, ..., n, \ n \ge 2,$ be
Q-independent random variables with values in $X$ and with
distributions $\mu_j$. Consider the linear forms $L_1=a_1\xi_1+
\cdots +a_n\xi_n$ and $L_2=b_1\xi_1 + \cdots + b_n\xi_n$, where the
coefficients $a_j, \ b_j$ are admissible integers for $X$ such that
$b_i a_j+ b_j a_i$ are admissible integers for $X$ for all $i, j$.
Assume that the conditional distribution of $L_2$ given $L_1$ is
symmetric. Then the following statements hold:}

\qquad \textit{$(i)$ If $X=\mathbb{R}^n\times D$, where $n\geq 0$
and $D$ is a torsion-free group, then all $\mu_j\in\Gamma(X)$;}

\qquad \textit{$(ii)$ If $X=X_{(3)}$ then all $\mu_j\in D(X)$.}

(II) \textit{If $X$ is not isomorphic to any of the groups mentioned
in (I), then there exist Q-independent random variables $\xi_j, j=1,
2, ..., n, n \ge 2,$ with values in $X$ and distributions $\mu_j$,
and admissible integers $a_j, \ b_j $ such that $b_i a_j + b_j a_i$
are admissible integers for $X$ for all $i, j$ such that the
conditional distribution of $L_2=b_1\xi_1 + \cdots + b_n\xi_n$ given
$L_1=a_1\xi_1+ \cdots +a_n\xi_n$ is symmetric, but all
$\mu_j\not\in\Gamma(X)*I(X)$.}

\medskip

\textbf{Proof.} {(I)}(i) If $X= \mathbb{R}^n\times D$, then
$Y\approx \mathbb{R}^n\times K$, where $K$ is a connected compact
group.  We can assume without loss of generality that
$Y=\mathbb{R}^n\times K$.

Lemma 7 implies that the characteristic functions of distributions
$\mu_j$ satisfy equation (\ref{1.1q}). Reasoning as in the beginning
of the proof of Theorem 1, we can assume that $\widehat\mu_j(y)> 0$
and the polynomial $q(u,v)$ is real-valued.

Consider the restriction of equation (\ref{1.1q}) to $K$. Since each
continuous polynomial is equal to a constant on compact elements,
the restriction of equation (\ref{1.1q}) to $K$ coincides with
equation (\ref{1.1}). It follows from Lemma 4 that
$\widehat\mu_j(y)=1$ on $K$. Then $\sigma(\mu_j)\subset A(X,K)=
\mathbb{R}^n$. Thus, we reduced the proof to the consideration of
equation (\ref{1.1q}) on $\mathbb{R}^n$. It follows from
\cite{My2018} that all $\mu_j\in\Gamma(X)$.

(ii) If $X=X_{(3)}$ then the connected component of zero of $X$ is
equal to zero. The independence and the Q-independence are coincides
on such groups. Therefore the statement of Theorem 4 in this case
follows from Theorem 1.

\textbf{(II)} Since independent random variables are Q-independent,
this part of Theorem 4 follows from part (II) of Theorem 2.
 $\blacksquare$

\end{document}